\documentclass[12pt]{article}

\usepackage{amsmath}
\usepackage{amssymb}

\newcommand{\eop}{\bigstar}  

\newcommand{\card}[1]{{\vert #1 \vert} }

\newcommand{\otp}[1]{\hbox{otp($#1$)}}

\newcommand{\cf}{{\rm cf}}

\newenvironment{Proof}{\noindent{\bf Proof.}}{\par\bigskip} 

\newtheorem{THEOREM}{Theorem}[section]

\newtheorem{Conclusion}[THEOREM]{Conclusion}

\newtheorem{Hypothesis}[THEOREM]{Hypothesis}

\newtheorem{LEMMA}[THEOREM]{Lemma}

\newtheorem{Main Theorem}[THEOREM]{Main Theorem}
\newenvironment{main Theorem}{\begin{Main Theorem}} 
{\end{Main Theorem}}

\newtheorem{Theorem}[THEOREM]{Theorem}

\newtheorem{Definition}[THEOREM]{Definition}

\newtheorem{Conventions}[THEOREM]{Conventions}

\newtheorem{Main Definition}[THEOREM]{Main Definition}
\newenvironment{main definition}{\begin{Main Definition}}
{\end{Main Definition}}

\newtheorem{Lemma}[THEOREM]{Lemma}

\newtheorem{Notation}[THEOREM]{Notation}

\newtheorem{Convention}[THEOREM]{Convention}

\newtheorem{Note}[THEOREM]{Note}

\newtheorem{Observation}[THEOREM]{Observation}

\newtheorem{Remark}[THEOREM]{Remark}

\newtheorem{Main Fact}[THEOREM]{Main Fact}
\newenvironment{main Fact}{\begin{Main Fact}}{\end{Main Fact}}

\newtheorem{Fact}[THEOREM]{Fact}

\newtheorem{Subfact}[THEOREM]{Subfact}

\newtheorem{Claim}[THEOREM]{Claim}

\newtheorem{Main Claim}[THEOREM]{Main Claim}
\newenvironment{main claim}{\begin{Main Claim}}{\end{Main Claim}}

\newtheorem{Crucial Claim}[THEOREM]{Crucial Claim}
\newenvironment{crucial claim}{\begin{Crucial Claim}}{\end{Crucial Claim}}

\newtheorem{Corrolary}[THEOREM]{Corrolary}

\newtheorem{Subclaim}[THEOREM]{Subclaim}

\newtheorem{Corollary}[THEOREM]{Corollary}

\newtheorem{Example}[THEOREM]{Example}

\newtheorem{Proposition}[THEOREM]{Proposition}

\newtheorem{Discussion}[THEOREM]{Discussion}

\newenvironment{Proof of the Subfact}
{\noindent{\bf Proof of the Subfact.}}{\par\bigskip}

\newenvironment{Proof of the Theorem}
{\noindent{\bf Proof of the Theorem.}}{\par\bigskip}

\newenvironment{Proof of the Conclusion}
{\noindent{\bf Proof of the Conclusion.}}{\par\bigskip}

\newenvironment{Proof of the Observation}
{\noindent{\bf Proof of the Observation.}}{\par\bigskip}

\newenvironment{Proof of the Fact}
{\noindent{\bf Proof of the Fact.}}{\par\bigskip}

\newenvironment{Proof of the Lemma}
{\noindent{\bf Proof of the Lemma.}}{\par\bigskip}

\newenvironment{Proof of the Claim}
{\noindent{\bf Proof of the Claim.}}{\par\bigskip}

\newenvironment{Proof of the Subclaim}
{\noindent{\bf Proof of the Subclaim.}}{\par\medskip}

\newenvironment{Proof of the Main Claim}
{\noindent{\bf Proof of the Main Claim.}}{\par\bigskip}

\newenvironment{Proof of the Crucial Claim}
{\noindent{\bf Proof of the Crucial Claim.}}{\par\bigskip}

\catcode`\@=11
\def\@begintheorem#1#2{\rm \trivlist \item[\hskip \labelsep{\bf #1\ #2.}]}
\def\@opargbegintheorem#1#2#3{\rm \trivlist
      \item[\hskip \labelsep{\bf #1\ #2\ (#3).}]}
\catcode`\@=12

\newcommand{\into}{\rightarrow}

\newcommand{\rest}{\upharpoonright}  

\newcommand{\cl}{\mathop{\rm cl}}
\newcommand{\cov}{\mathop{\rm cov}}

\newcommand{\ran}{\mathop{\rm ran}}

\newcommand{\deq}{\buildrel{\rm def}\over =}

\newcommand{\PP}{{\cal P}}
\newcommand{\RR}{{\cal R}}

\newcount\skewfactor
\def\mathunderaccent#1#2 {\let\theaccent#1\skewfactor#2
\mathpalette\putaccentunder}
\def\putaccentunder#1#2{\oalign{$#1#2$\crcr\hidewidth
\vbox to.2ex{\hbox{$#1\skew\skewfactor\theaccent{}$}\vss}\hidewidth}}

\title{On club-like principles on regular cardinals above $\beth_\omega$}

\author{Mirna D\v zamonja\\
School of Mathematics\\
University of East Anglia\\
Norwich, NR4 7TJ\\
UK\\
\scriptsize{h020@uea.ac.uk}\\
\scriptsize{http://www.mth.uea.ac.uk/people/md.html}}

\begin{document}
\maketitle
\begin{abstract} We prove that for regular $\lambda$
above a strong limit singular $\mu$ certain guessing principles follow just from
cardinal arithmetic assumptions. The main result is that for such $\lambda$ and
$\mu$
there are coboundedly many regular $\kappa<\mu$ such that
$\clubsuit^-(S^\lambda_\kappa)$ holds whenever $\lambda=\lambda^{<\kappa}$.
{\footnote{The author thanks EPSRC for their support through the grant number
GR/M71121. She is also grateful to Prof. Balcar and Prof. Simon of the Charles
University in Prague for their invitation to speak at the 29th Winter School of
Abstract Analysis, which proved to be a most productive and pleasant experience.}}
\end{abstract}

\section{Introduction.}
The main result of this note is that for any regular cardinal $\lambda$ above
$\beth_\omega$ there are unboundedly many regular cardinals $\kappa<\beth_\omega$ such that provided
$\lambda=\lambda^{<\kappa}$, the combinatorial principle
$\clubsuit^-(S^\lambda_\kappa)$ holds. That principle is defined in the
following Definition \ref{club}, and the notation $S^\lambda_\kappa$ is
recalled in \ref{oznake}.

\begin{Definition}\label{club} Suppose that $S$ is a stationary subset of a
regular cardinal $\lambda$. Then $\clubsuit^-(S)$ is the statement claiming the
existence of
a sequence $\langle \PP_\delta:\,\delta\in S\rangle$
such that 
\begin{description}
\item{(i)} each $\PP_\delta$ is a family of $<\lambda$ many subsets of
$\delta$ and 
\item{(ii)} for every unbounded subset $A$ of $\lambda$
there are stationarily many $\delta$ such that for some $X\in \PP_\delta$ with
$\sup(X)=\delta$ we have
$X\subseteq A$. 
\end{description}
\end{Definition}

We also prove that some other
similar
combinatorial principles on such $\lambda$ follow just from the assumptions on
cardinal arithmetic. In fact, the same theorems hold more generally in a
situation in which $\beth_\omega$ is replaced by any strong limit singular
cardinal.
Our proofs are an application of a (consequence of) a powerful theorem of
Shelah in \cite{Sh 460}, Theorem \ref{460} below. 
The methods are similar to the ones used in \cite{Sh 460}
to prove e.g. that for $\lambda$ as above the assumption $\lambda^{<\lambda}=\lambda$
implies that $\diamondsuit^-$ holds at $\lambda$. 
 
Throughout the note we use the notation given below. Note that $\cov$ as used
here is a special case of a more general notation used in pcf theory, but to increase
readability we only quote the instance of it that we actually use.

\begin{Notation}\label{oznake} Suppose that $\kappa$ is a regular cardinal and
$\alpha>\kappa$ an ordinal. Then
\begin{description}
\item{(1)}
$S^\alpha_\kappa\deq\{\beta<\alpha:\,\cf(\beta)=\kappa\}$,
\item{(2)}
$S^\alpha_{<\kappa}\deq\{\beta<\alpha:\,\cf(\beta)<\kappa\}$,
\item{(3)}
\begin{multline*}
\cov(\alpha,\kappa^+,\kappa^+,\kappa)\deq
\min\{\theta:\,(\exists\PP
\subseteq [\alpha]^{\le\kappa})\\
\card{\PP}=\theta\,\,\&\,\,(\forall A\in
[\alpha]^{\le\kappa})(\exists X\in [\PP]^{<\kappa}) A\subseteq \bigcup X\}.
\end{multline*} 
\item{(4)} For a subset $A$ of $\kappa$ we let
$\lim(A)\deq\{\beta<\kappa:\,\beta=\sup(A\cap\beta)\}$ and $\cl(A)=A\cup
\lim(A)$.
\end{description}
\end{Notation}

The theorem we need for our application is given below as 
Theorem \ref{460}. Its statement is modulo the notation an easy consequence of
Theorem 1.1. of \cite{Sh 460} combined with another deep theorem of cardinal
arithmetic ,
the `cov versus pp' theorem of Shelah. As this may not be immediate from
reading
\cite{Sh 460}, for the benefit of an interested reader
we briefly comment on how the connection can be seen. 

\begin{Theorem}\label{460} (Shelah) Suppose that $\mu$ is a strong limit
singular cardinal. Then
for $\lambda>\mu$, for every regular large enough $\kappa<\mu$, we
have that for all $\alpha <\lambda$,
\[
\mbox{cov}(\alpha,\kappa^+,\kappa^+,\kappa)<\lambda.
\]
\end{Theorem}

{\em Sketch of the proof}. The statement of Theorem 1.1. of \cite{Sh 460} is
that in the situation as described by the assumptions of Theorem \ref{460},
there are only boundedly many $\kappa<\mu$ such that for some $\lambda^\ast\in
(\mu,\lambda)$ we have ${\rm pp}_{\Gamma(\mu^+, \kappa)}(\lambda^\ast)\ge\lambda$.
The notation to the extent needed here will be described below.

Suppose $\alpha<\lambda$. As clearly
$\mbox{cov}(\alpha,\kappa^+,\kappa^+,\kappa)=\mbox{cov}(\card{\alpha},
\kappa^+,\kappa^+,\kappa)$ for any $\kappa$, we can assume that $\alpha$ is a
cardinal $\theta$. Let $\kappa<\mu$ be large enough uncountable such that for
no $\lambda^\ast\in (\mu,\lambda)$ do we have
${\rm pp}_{\Gamma(\mu^+, \kappa)}(\lambda^\ast)\ge\lambda$. The notation used
here is that for a cardinal $\sigma$
\begin{multline*}
\Gamma(\sigma,\kappa)\deq\{I:\,\mbox{ for some cardinal }\theta_I<\sigma\\
I\mbox{ is a proper }\kappa\mbox{-complete ideal on }\theta_I\}
\end{multline*}
and
\begin{multline*}
{\rm pp}_{\Gamma(\mu^+, \kappa)}(\lambda^\ast)
\deq\sup\{\mbox{tcf}(\Pi\frak a/J):\,\frak a\mbox{ is a set of regular 
cardinals}\\
\mbox{unbounded in }\lambda^\ast,\\
J\in \Gamma(\mu^+,\kappa)\mbox{ and }\mbox{tcf}(\Pi\frak a/J)
\mbox{ is well defined}\}.
\end{multline*}
For our purposes here it is not important what the notation
$\mbox{tcf}(\Pi\frak a/J)$ means exactly, one should simply observe that
$\Gamma(\kappa^+,\kappa)\subseteq\Gamma(\mu^+,\kappa)$
and hence ${\rm pp}_{\Gamma(\mu^+, \kappa)}(\lambda^\ast)\ge {\rm
pp}_{\Gamma(\kappa^+, \kappa)}(\lambda^\ast)$ for all relevant
$\lambda^\ast$.
This implies that
for no $\lambda^\ast\in (\mu,\lambda)$ do we have ${\rm pp}_{\Gamma(\kappa^+, \kappa)}
(\lambda^\ast)\ge\lambda$.

Now we quote Shelah's `cov versus pp' theorem, \cite {Sh g}, II 5.4., which
says that
\[
\mbox{cov}(\theta,\kappa^+,\kappa^+,\kappa)+\theta=\sup
\{\mbox{pp}_{\Gamma(\kappa^+,\kappa)}(\lambda^\ast):\,\lambda^\ast\in [\kappa,
\theta]\}+\theta,
\]
leading us to conclude that
$\mbox{cov}(\theta,\kappa^+,\kappa^+,\kappa)<\lambda$.
$\eop_{\ref{460}}$

\medskip

We shall also use another staple of cardinal arithmetic, namely the club
guessing principle quoted in the following

\begin{Theorem}\label{clubg} (Shelah, \cite{Sh g}, III,\S2) Suppose that
$\aleph_0<\cf(\kappa)=\kappa$ and $\kappa^+<\lambda=\cf(\lambda)$. Then there
is a sequence $\bar{e}=\langle e_\delta:\,\delta\in S^\lambda_\kappa\rangle$ of
sets such that for each $\delta$ we have $\otp{e_\delta}=\kappa$ and $e_\delta$ is a
club subset of $\delta$ consisting of points of cofinality $<\kappa$, and for
every club $E$ of $\lambda$ there are stationarily many $\delta$ such that
$e_\delta\subseteq E$.

If $\kappa=\aleph_0$, then there is a sequence $\bar{e}$ of the above form such
that each $e_\delta$ is a cofinal subset of $\delta$ of order type $\omega$,
and for
every club $E$ of $\lambda$ there are stationarily many $\delta$ such that
$e_\delta\subseteq E$.
\end{Theorem}

\section{The results.}\label{clubs}
To simplify the notation, which involves dealing with many cardinals at a time,
we first formulate and prove the theorem in lesser generality where
$\beth_\omega$ is the strong limit singular we work with. The same proof gives
the fully general result, as indicated in Theorem \ref{main}.

\begin{Theorem}\label{guessing} Suppose that $\lambda$ is a regular cardinal
with $\lambda>\beth_\omega$. 

Then there are coboundedly many regular $\kappa<\beth_\omega$ such that the
following statements hold:

\begin{description}
\item{(1)} If $\lambda^{<\kappa}=\lambda$, then
$\clubsuit^-(S^\lambda_\kappa)$ holds. Precisely, 
there is a sequence $\langle\PP_\delta:\,\delta\in S^\lambda_\kappa\rangle$ such that
\begin{description}
\item{(i)} each $\PP_\delta$ is a family of $<\lambda$ elements of
$[\delta]^{\le\kappa}$ and
\item{(ii)} for every $A\in [\lambda]^\lambda$, there are stationary many
$\delta$ such that for some $X$ in $\PP_\delta$ we have $X\subseteq A$
and $\sup(X)=\delta$.
\end{description}

\item{(2)} There is a sequence $\langle\PP^0_\delta:\,\delta\in
S^\lambda_\kappa\rangle$ satisfying (1)(i) above and such that for all $A\in
[\lambda]^\lambda$ there is a club $E$ of $\lambda$ such that for every
$\delta\in E\cap S^\lambda_\kappa$, for some $a\in \PP^0_\delta$ we have
$\sup(A\cap a)=\delta$.

\item{(3)} If $\theta<\lambda\implies \theta^{<\kappa}<\lambda$, then there is
a sequence $\langle\RR_\delta:\,\delta\in
S^\lambda_\kappa\rangle$ satisfying (1)(i) above and
\begin{description}
\item{(ii)${}^+$} for every sequence $\langle
a_\delta:\,\delta\in S^\lambda_\kappa\rangle$ of sets such that each $a_\delta$
is a subset of $\delta$ of order type $\le\kappa$, there is a club $C$
of $\lambda$ such that $\delta\in C\cap S^\lambda_\kappa\implies a_\delta\in
\RR_\delta$. \end{description}
\end{description}
\end{Theorem}

\begin{Proof} For $\alpha<\lambda$, let
\[
R_\alpha\deq\{\kappa\mbox{ regular
}<\beth_\omega:\,\cov(\alpha,\kappa^+,\kappa^+,\kappa)<\lambda\}.
\]
By Theorem \ref{460}, for each such $\alpha$ there is
$n_\alpha<\omega$ such that $R_\alpha$ contains all regular cardinals in the
interval $[\beth_{n_\alpha},\beth_\omega)$. Hence there is $n^\ast<\omega$
such that for unboundedly many $\alpha<\lambda$ we have that $n_\alpha=n^\ast$.
As it is easily seen that 
\[
\alpha<\beta\implies \cov(\alpha,\kappa^+,\kappa^+,\kappa)\le
\cov(\beta,\kappa^+,\kappa^+,\kappa),
\]
it follows that for all $\alpha<\lambda$, the set $R_\alpha$ contains
all regular cardinals in $[\beth_{n^\ast},\beth_\omega)$. Let us fix a regular
cardinal $\kappa>\aleph_0$ in the interval $[\beth_{n^\ast},\beth_\omega)$ and
show that all three statements of the Theorem hold for such $\kappa$.

For each $\alpha<\lambda$ let $\PP_\alpha^0$ be a family exemplifying
that $\cov(\alpha,\kappa^+,\kappa^+,\kappa)<\lambda$. 
The sequence needed for (2) is in fact $\langle \PP^0_\delta:\,\delta\in S^\lambda_\kappa
\rangle$, a point to which we shall briefly return later,
but for the moment we go on to the main part of the proof, which is the proof
of (1).

{\em Proof of (1)}. As we are assuming 
$\lambda^{<\kappa}=\lambda$, let us enumerate $[\lambda]^{<\kappa}=
\{A^\ast_i:\,i<\lambda\}$ so that each set in the enumeration appears
$\lambda$ many times. For $\delta\in S^\lambda_\kappa$ let 
\[
\PP^1_\delta=\{(\bigcup_{i\in B}A^\ast_i)\cap\delta:\,B\in\PP^0_\delta\},
\]
hence each $X\in \PP^1_\delta$ is a subset of $\delta$ of size $\le\kappa$
and $\card{\PP^1_\delta}<\lambda$. Fixing $\delta\in S^\lambda_\kappa$ for a
moment, we have that for each $X\in \PP^1_\delta$ the size of $X$ is
$\le\kappa$, so the size of $\PP(X)$ is $\le 2^\kappa<\beth_\omega<\lambda$,
leading us to conclude that 
\[
\PP_\delta\deq\{Y:\,(\exists X\in \PP^1_\delta) Y\subseteq X\}
\]
also has size $<\lambda$. We shall proceed to show that
$\langle \PP_\delta:\,\delta\in S^\lambda_\kappa\rangle$ is a sequence
as required. Part (i) of our requirement is clearly satisfied, so let us
proceed to part (ii). For this we shall first need to fix a club guessing 
sequence $\langle e_\delta:\,\delta\in S^\lambda_\kappa\rangle $ as provided by
Theorem \ref{clubg}. For each $\delta\in S^\lambda_\kappa$, let $e_\delta=
\{\zeta^\delta_\gamma:\,\gamma<\kappa\}$ be the increasing enumeration of
$e_\delta$.

Let $A\in [\lambda]^{\lambda}$ be given. For $\varepsilon \in
S^\lambda_{<\kappa}$ define $X_\varepsilon=X_\varepsilon^A$ to be a subset of
$A$ of size $<\kappa$ with $\sup(X_\varepsilon)=\varepsilon$, if such a set
exists. Now define a function
$h_A:\,S^\lambda_{<\kappa}\into\lambda$ by the following recursive definition
\begin{equation*}
h_A(\varepsilon)\deq
\begin{cases}
\min\{i>\sup_{\beta\in
S^\varepsilon_{<\kappa}}h_A(\beta):\,X_\varepsilon=A^\ast_i\}
&
\text{ if $X_\varepsilon$ is defined},\\
\sup_{\beta\in
S^\varepsilon_{<\kappa}}h_A(\beta)+1
&
\text{ otherwise}.
\end{cases}
\end{equation*}
Let $E\deq\lim(\{\delta<\lambda:\,(\forall\varepsilon\in S^\delta_{<\kappa})
h_A (\varepsilon)\le\delta\})$, hence a club of $\lambda$. Note that if
$\delta\in E\cap S^\lambda_\kappa$, then for all $\varepsilon\in
S^\delta_{<\kappa}$
we actually have $h_A (\varepsilon)<\delta$. Let us choose $\delta\in E\cap
S^\lambda_\kappa$
such that $e_\delta\subseteq\lim(A)$. This in particular means that for every
$\gamma<\kappa$ the set $X_{\zeta^\delta_\gamma}$ has been
defined. For such $\gamma$, let $i_\gamma\deq
h_A(\zeta^\delta_\gamma)$, hence $\langle i_\gamma:\,\gamma<\kappa\rangle$ is
a strictly increasing sequence and for each $\gamma$ we have
$A^\ast_{i_\gamma}=X_{\zeta^\delta_\gamma}$. As $\{i_\gamma:\,\gamma<\kappa\}
\in [\delta]^\kappa$, there are sets $\{B_j:\,j<j^\ast<\kappa\}$ in
$\PP_\delta^0$ such that
$\{i_\gamma:\,\gamma<\kappa\}\subseteq\bigcup_{j<j^\ast} B_j$. By the
regularity of $\kappa$, there is $j<j^\ast$ such that
$\card{\{i_\gamma:\,\gamma<\kappa\}\cap B_j}=\kappa$. Let $B=B_j$ for some such
$j$.

Consider $(\bigcup_{i\in B}A^\ast_i)\cap\delta$. Clearly, this set is a
superset of $\bigcup_{i_\gamma\in B}X_{\zeta^\delta_\gamma}$ (so it has size
$\kappa$) and is a member of $\PP^1_\delta$. For this reason,
$\bigcup_{i_\gamma\in B}X_{\zeta^\delta_\gamma}\in \PP_\delta$, and this set is
clearly an unbounded subset of $A\cap\delta$ of size $\kappa$.

{\em Proof of (2)}. This follows trivially with $\langle
\PP_\delta^0:\,\delta\in S^\lambda_\kappa\rangle$ as above, since
by taking
$\delta\in \lim(A)\cap S^\lambda_\kappa$, we obtain that $A\cap\delta$ is unbounded in $\delta$ and
covered by $<\kappa$ many elements of $\PP^0_\delta$. Hence, by the regularity
of $\kappa$ we obtain that there is an element $X$ of $\PP^0_\delta$ with
$\sup(A\cap X)=\delta$.

{\em Proof of (3)}. For each relevant $\delta$, we form the family
$\PP^0_\delta$ as in the proof of (1). Fixing
$\delta\in S^\lambda_\kappa$ for a moment and letting
$\theta=\card{\PP^0_\delta}$, we have $\theta^{<\kappa}<\lambda$, so we can let
$\PP^2_\delta$ consist
of the unions of all subfamilies of $\PP^1_\delta$ which have size $<\kappa$
and obtain a family of elements of elements of $[\delta]^{\le\kappa}$ of size
$<\lambda$. The proof now follows the proof of (1), but we give the details for
the sake of completeness.

As $\theta<\lambda\implies\theta^{<\kappa}<\lambda$ and $\lambda$ is regular,
we have
$\lambda^{<\kappa}=\lambda$. We enumerate
$[\lambda]^{<\kappa}=\{A^\ast_i:\,i<\lambda\}$. Let
$\PP^3_\delta\deq\{(\bigcup_{i\in B}A^\ast_i)\cap\delta:\,B\in \PP^2_\delta\}$,
and let $\RR_\delta\deq\{Y:\,(\exists X\in \PP^3_\delta)X\subseteq Y\}$, for
each relevant $\delta$.
Let 
\[
E\deq\cl(\{\beta<\lambda:\,(\forall X\in [\beta]^{<\kappa})(\exists i<\beta)
X=A^\ast_i\}).
\]
Note that if $\delta\in E\cap S^\lambda_\kappa$ then for all $X\in
[\delta]^{<\kappa}$ we have $X=A^\ast_i$ for some $i<\delta$. We claim that for
each such $\delta$ the set $a_\delta$ is in $\PP_\delta$. Let
$f_\delta:\,\kappa\into a_\delta$ be the increasing enumeration of $a_\delta$
and for $\gamma<\kappa$ let $X_\gamma=\ran(f_\delta\rest\gamma)$. For each such
$\gamma$ let $i_\gamma<\delta$ be such that $X_\gamma=A^\ast_{i_\gamma}$. Hence
there are sets $\{B^\ast_j:\,j<j^\ast<\kappa\}$ in $\PP^0_\delta$ such that
$\{i_\gamma:\,\gamma<\kappa\}\subseteq\bigcup_{j<j^\ast}B_j\deq B$. We have
that $B\in \PP^2_\delta$, hence $(\bigcup_{i\in B}A^\ast_i)\cap\delta\in
\PP^3_\delta$ and is a superset of $a_\delta$, so $a_\delta\in \RR_\delta$.
$\eop_{\ref{guessing}}$
\end{Proof}

A more general theorem is 

\begin{Theorem}\label{main} The analogue of Theorem \ref{guessing} holds when
$\beth_\omega$ is replaced by any other strong limit singular cardinal $\mu$.
\end{Theorem}

\begin{Proof} Exactly the same as that of Theorem \ref{guessing}, using the
full generality of Theorem \ref{460} and replacing $\beth_\omega$ by $\mu$
throughout.
$\eop_{\ref{main}}$
\end{Proof}

\section{Concluding Remarks.} The main result we proved is that when $\mu$ is a
strong limit singular cardinal and $\lambda$ is a regular cardinal above $\mu$,
there are coboundedly many regular $\kappa<\mu$ such that 
\[
\lambda=\lambda^{<\kappa}\implies\clubsuit^{-}(S^\lambda_\kappa),
\]
hence the existence of the guessing sequence follows simply from the cardinal
arithmetic assumed. When combined with the result of Shelah in \cite{Sh 460}
which under these conditions shows the equivalence of the assumption
$\lambda=\lambda^{<\lambda}$ with $\diamondsuit^-$, an immediate consequence is
that $\clubsuit^{-}(\lambda)$ and $\diamondsuit^-(\lambda)$ are different,
a fact whose analogue at $\omega_1$ requires a rather serious proof (Shelah,
see \cite{Sh f} e.g).
In fact our result implies the former among the Shelah's results, as it
is well known that $\clubsuit^-(\lambda)+\lambda^{<\lambda}=\lambda
\implies \diamondsuit^-(\lambda)$.
It would be interesting to know if when we in addition
assume that $\lambda$ as above is a successor cardinal, then 
$\lambda=\lambda^{<\kappa}\implies\clubsuit$. The analogue of this for
$\diamondsuit$ follows from the above mentioned result of \cite{Sh 460}
and Kunen's argument on the equivalence between $\diamondsuit$ and
$\diamondsuit^-$
at successor cardinals (see \cite{Kune} e.g). We have the impression that the
answer to the question is negative, since it is known
by \cite{DzSh} that $\clubsuit^-$ and $\clubsuit$ differ at $\aleph_1$.

\eject

\end{document}